\documentclass{amsart}
\usepackage{amssymb,latexsym}

\numberwithin{equation}{section}

\newtheorem{theorem}{Theorem}
\newcommand{\bt}{\begin{theorem}}
\newcommand{\et}{\end{theorem}}
\newtheorem{lemma}{Lemma}
\newcommand{\bl}{\begin{lemma}}
\newcommand{\el}{\end{lemma}}
\newtheorem{corollary}{Corollary}
\newcommand{\bc}{\begin{corollary}}
\newcommand{\ec}{\end{corollary}}
\newtheorem{problem}{Problem} 
\newcommand{\bprob}{\begin{problem}}
\newcommand{\eprob}{\end{problem}}
\newcommand{\beq}{\begin{equation}}
\newcommand{\eeq}{\end{equation}}
\newcommand{\benum}{\begin{enumerate}}
\newcommand{\eenum}{\end{enumerate}}

\newcommand{\PP}{\ensuremath{\mathbf P }}
\newcommand{\N}{\ensuremath{ \mathbf N }}
\newcommand{\Q}{\ensuremath{ \mathbf{Q} }}
\newcommand{\Z}{\ensuremath{\mathbf Z}}

\begin{document}

\title[Phase transitions in infinitely generated groups]{Phase transitions in infinitely generated groups, and related problems in additive number theory}
\author{Melvyn B. Nathanson}
\address{Department of Mathematics, CUNY (Lehman College and the Graduate Center)}
\email{melvyn.nathanson@lehman.cuny.edu}
\curraddr{Department of Mathematics, Princeton University}
\email{melvyn@princeton.edu}


\date{}

\begin{abstract}
Let $A$ be an infinite set of generators for a group $G$, and let $L_A(r)$ denote the number of elements of $G$ whose word length with respect to $A$ is exactly $r$.  The growth function $L_A$ is a function from the nonnegative integers $\N_0$ to the set $\N_0 \cup \{\infty\}$.  The purpose of this note is to determine all growth functions associated to infinite generating sets for groups, and to describe a phase transition phenomenon associated with infinite generating sets.  A list of open problems is also.included.  
\end{abstract}

\maketitle 

\section{Generating sets and the word growth function}
Let $\N = \{1,2,3,\ldots\}$ denote the set of positive integers, $\N_0 = \N \cup \{ 0\}$  the set of nonnegative integers, and \Z\ the set of all integers.  Let $\N_{\infty} = \N \cup \{\infty\}$.  
We denote the cardinality of the set $X$ by $|X|$.

Let $G$ be a group, written multiplicatively, and let $A$ be a finite or infinite subset of $A$.  Let $A^{-1} = \{ a^{-1} : a \in A \}$.   The set $A$ is \emph{symmetric} if $A = A^{-1}$.  The \emph{subgroup generated by $A$}, denoted $\langle A \rangle$, is the set of elements of $G$ that can be written as a finite product of elements of $A$ and their inverses, that is, as a finite product of elements of $A \cup A^{-1}$.  We define the empty product to be the  identity $e$ of $G$.   The set $A$ is  called a \emph{set of generators for the group $G$} if $G = \langle A \rangle$.  In this paper we study infinite generating sets.  Of course, if the group $G$ happens to be finitely generated, then every infinite generating set contains a finite subset that also generates $G$, but we do not require that our generating sets be minimal.

Let $A$ be a set of generators for the group $G$.  For every $x \in G$, the \emph{word length of $x$ with respect to $A$}, denoted  $\ell_A(x)$, is the smallest integer $r$ such that $x$ can be represented as a product of $r$ elements of $A \cup A^{-1}$. 
For $x \in G$, we have $\ell_A(x) = 0$ if and only if $x = e$.  Also, $\ell_A(x) = 1$ if and only if $x \neq e$ and $x \in A \cup A^{-1}$. 
The word length function is symmetric in the sense that 
\beq   \label{PhaseTransition:Symmetry}
\ell_A(x) = \ell_A(x^{-1})
\eeq
for all $x \in G$, and satisfies the triangle inequality
\beq   \label{PhaseTransition:TriangleInequality}
\ell_A(xy) \leq \ell_A(x) + \ell_A(y)
\eeq
for all $x,y \in G$.  For every nonnegative integer $r$, we define the \emph{sphere of radius $r$}
\[
S_A(r) = \{ x \in G : \ell_A(x) = r \}
\]
and the \emph{word growth function} 
\[
L_A(r) = \left| S_A(r) \right|.
\]

 If $G$ is an additive group with generating set $A$, then we denote by $-A$ the set of inverses of elements of $A$, and we define the spheres $S_A(r)$ and the functions $\ell_A(x)$ and $L_A(r)$ analogously.

Here are two examples of growth functions.  Let $\Q^+$ denote the multiplicative group of positive rational numbers.  The set $\PP$ of prime numbers generates $\Q^+$.  The sphere $S_{\PP}(r)$ consists of all positive rational numbers that can be written in the form $p_1^{k_1}p_2^{k_2}\cdots p_{\ell}^{k_{\ell}}$, where $p_1,p_2,\ldots, p_{\ell}$ are distinct primes, $k_1,k_2,\ldots, k_{\ell}$ are nonzero integers, and $\sum_{i=1}^{\ell} |k_i| = r$.   It follows that $L_{\PP}(r) = \infty$ for all $r\geq 1$.  

Let \Z\ be the additive group of integers, and let $m$ be a positive integer.   The set $A = \{1\} \cup \{m, 2m, 3m, 4m,\ldots\}$ is an infinite generating set for \Z, and  
\[
L_A(r) = 
\begin{cases}
\infty & \text{if $1 \leq r \leq \left[\frac{m}{2}\right]  +1 $ } \\
 0 & \text{if $r > \left[\frac{m}{2}\right]  + 1 $ }. 
\end{cases}
\]

Let $A$ be an infinite generating set for a group $G$.  We shall prove that either $L_A(r) = \infty$ for all $r \geq 1$, or there exist numbers $r \in \N$ and  $s \in \N_{\infty}$  such that $L_A(r') = \infty$ for $1 \leq r' < r$, $L_A(r) = s$, and  $L_A(r') = 0$ for $r' > r$.  The ordered pair $(r,s)$ is called the \emph{phase transition} associated to the group $G$ and generating set $A$, and the sphere $S_A(r)$ is called the \emph{transition set}.  
This is the last nonempty set in the sequence $\{S_A(r')\}_{r'=0}^{\infty}$.  Note that $r=1$ if and only if $s=\infty$ and $G = A \cup A^{-1} \cup \{e\}$.    Also, $r=2$ if and only if $S_A(2) = G \setminus \left( A \cup A^{-1} \cup \{e\} \right) \neq \emptyset$.

The phase transition $(r,s)$  is \emph{finite} if $s\in \N$ and \emph{infinite} if $s = \infty$.  We construct examples  to show that every pair of numbers $(r,s)$ with $r \in \N$, $r \geq 2$,  and  $s \in \N_{\infty}$ can occur as the phase transition associated to a generating set for a group.  Moreover, for the additive group \Z\ of integers, we prove that for every integer $r \geq 2$ and for  every finite symmetric set $S$ of nonzero integers with $|S|=s$, there is an infinite generating set $A$ such that, with respect to $A$, the set $S$ is the transition set for a phase transition of the form $(r,s)$.

\section{Existence of phase transitions}
In this section we prove the fundamental theorem on the existence of phase transitions for groups with an infinite set of generators.

\bl    \label{PhaseTransition:lemma:translate}
Let $A$ be a generating set for a group $G$.  
If $x , y\in G$, then 
\beq    \label{PhaseTransition:translate}
\ell_A(x) - \ell_A(y) \leq \ell_A(xy) \leq \ell_A(x) + \ell_A(y).
\eeq
If $x \in G$ and $a\in A$, then 
\beq  \label{PhaseTransition:translateLeft}
\left| \ell_A(ax) - \ell_A(x) \right| \leq 1
\eeq
and 
\beq  \label{PhaseTransition:translateRight}
\left| \ell_A(xa) - \ell_A(x) \right| \leq 1.
\eeq
\el

\begin{proof}
By the symmetry relation~\eqref{PhaseTransition:Symmetry} and  the triangle inequality~\eqref{PhaseTransition:TriangleInequality}, we have
\[
\ell_A(x) = \ell_A(xyy^{-1}) \leq \ell_A(xy) + \ell_A(y^{-1}) 
=  \ell_A(xy) + \ell_A(y)
\]
and so
\[
\ell_A(x) - \ell_A(y) \leq \ell_A(xy) \leq \ell_A(x) + \ell_A(y).
\]
This proves~\eqref{PhaseTransition:translate}.  

If $e=a\in A$, then $\ell_A(ax) = \ell_A(x) = \ell_A(xa)$ and inequalities~\eqref{PhaseTransition:translateLeft} 
and~\eqref{PhaseTransition:translateRight} hold.
If $a\in A\setminus \{ e\}$, then $\ell_A(a) = 1$ and inequality~\eqref{PhaseTransition:translate} with $y=a$ gives 
\[
-1 \leq  \ell_A(xa) -  \ell_A(x) \leq 1.
\]
Similarly,  inequality~\eqref{PhaseTransition:translate} with $x=a$ gives  
\[
-1 \leq  \ell_A(ay) -  \ell_A(y) \leq 1.
\]
This completes the proof.
\end{proof}

\bl   \label{PhaseTransition:lemma:u-sequence} 
Let $u_0, u_1,\ldots, u_k$ be a sequence of  integers such that $u_0 \leq u_k$ and $|u_i - u_{i-1}| \leq 1$ for $i=1,\ldots, k$.  If $r\in \Z$ and $u_0 \leq r \leq u_k$, then there exists a nonnegative  integer $j \leq k$ such that $u_j = r$.
\el

\begin{proof}
If $u_0 = r$, let $j=0$.  If $u_0 < r$, let $j$ be the greatest positive integer such that $u_i < r$ for $i=0,1,\ldots, j-1$.  We have $j \leq k$ since  $u_k \geq r$.  Since $u_{j-1} \leq r-1$ and $u_j \geq r$, the inequalities
\[
u_j - u_{j-1} = |u_j - u_{j-1}| \leq 1
\]
and
\[
u_{j-1} \leq r-1 < r \leq u_j \leq u_{j-1}+1
\]
imply that $u_j=r$.  
\end{proof}

\bl     \label{PhaseTransition:lemma:r-sequence} 
Let $A$ be a generating set for a group $G$.  
Let $a_i \in A$ for $i=1,\ldots, k$.   If the nonnegative integer $r$ satisfies $r \leq \ell_A(a_1a_2\cdots a_k)$, then there exists  a nonnegative  integer $j \leq k$ such that $r = \ell_A(a_1a_2\cdots a_j)$.
\el

\begin{proof}
Apply Lemmas~\ref{PhaseTransition:lemma:translate} and~\ref{PhaseTransition:lemma:u-sequence} with $u_i =  \ell_A(a_1a_2\cdots a_i)$ for $i=0,1,\ldots, k$.
\end{proof}

\bl   \label{PhaseTransition:lemma:cut}
Let $A$ be a generating set for a group $G$.  
Let $a_i \in A \cup A^{-1}$ for $i = 1,2,\ldots,r$.  
If
\[
\ell_A(a_1a_2\cdots a_r) = r
\]
and $1 \leq i \leq j \leq r$, then 
\[
\ell_A(a_ia_{i+1}\cdots a_{j}) = j-i+1.
\]
\el

\begin{proof}
Let $x=a_1a_2\cdots a_r$ and $x' = a_ia_{i+1}\cdots a_{j}$.  
We have
\[
\ell_A(x') = \ell_A(a_ia_{i+1}\cdots a_{j}) \leq j-i+1.
\]
Since
\[
x = a_1a_2\cdots a_{i-1} x' a_{j+1}\cdots a_r
\]
the triangle inequality implies that 
\begin{align*}
r & = \ell_A(x) = \ell_A( a_1a_2\cdots a_{i-1} x' a_{j+1}\cdots a_r ) \\
& \leq \ell_A( a_1a_2\cdots a_{i-1} ) + \ell_A(  x'  ) + \ell_A( a_{j+1}\cdots a_r )\\
& \leq (i-1) + \ell_A(x') + (r-j)
\end{align*} 
and so 
\[
\ell_A(x') \geq j-i+1.
\]
This completes the proof.
\end{proof}

\bt           \label{PhaseTransition:theorem:fundamental}
Let $A$ be an infinite generating set for a group $G$.  
If $r > 1$ and $L_A(r) < \infty$, then $L_A(r') = 0$ for all $r' > r$.  
\et

\begin{proof}
Recall that $S_A(1) = (A \cup A^{-1})\setminus\{e\}$, and so $L_A(1) = |S_A(1)|  = \infty$.  
We define the sets
\[
S_A(<r) = \{ x\in G : \ell_A(x) < r \}
\]
and
\[
S_A(>r) = \{ x\in G : \ell_A(x) > r \}.
\]
Then $S_A(1) \subseteq S_A(<r)$.  
We shall prove that $S_A(>r)$ is the empty set.  Suppose not.  If $x \in S_A(>r)$, then  $x$ is an element of $G$ whose length is  $r' > r$, and there exist elements $a_i \in (A \cup A^{-1})\setminus \{ e\}$ for $i = 1,\ldots, r'$ such that 
\[
x = a_1a_2\cdots a_{r'}.
\]
By Lemma~\ref{PhaseTransition:lemma:translate}, for every $a\in A$ the element 
\[
ax = aa_1a_2\cdots a_{r'}
\]
has length at least $r'-1\geq r$.  The generator $a$ has length 1.
Consider the sequence of elements
\[
a, aa_1, aa_1a_2, a a_1 a_2 a_3, \ldots, aa_1a_2 \cdots a_{r'}.
\]
The first term in this sequence has length 1 and the last term has length at least $r$.  By Lemma~\ref{PhaseTransition:lemma:r-sequence}, at least one term of this sequence has length exactly $r$.  Since $r > 1$, this term must be different from $a$.  It follows that there exists a function $f:A \rightarrow \{1,2,\ldots, r\}$ so that, if $f(a) = i$, then $aa_1\cdots a_i$ has length $r$.  Since the generating set $A$ is infinite,  the pigeonhole principle implies that there exists an integer $i$ such that $f^{-1}(i) = \{a\in A : f(a) = i\}$ is an infinite set.  If $a \in f^{-1}(i)$, then $aa_1\cdots a_i$ has length $r$.   If $a,a' \in f^{-1}(i)$ and $a \neq a'$,  then $aa_1\cdots a_i \neq a'a_1\cdots a_i$, and so the sphere $S_A(r)$ is infinite.  This is impossible, hence $S_A(>r)$ is empty.
\end{proof}

\bc       \label{PhaseTransition:corollary:fundamental}
Let $A$ be an infinite generating set for a group $G$.  If $r \in \N$ and $0 < L_A(r) < \infty$, then $r > 1$ and $L_A(r') = \infty$ for $r' < r$ and $L_A(r') = 0$ for $r' > r$.  If $L_A(r) = 0$ or $\infty$ for all positive integers $r$, then either $L_A(r) = \infty$ for all $r$, or there is a positive integer $r$ such that $L_A(r') = \infty$ for all positive integers $r' \leq r$ and $L_A(r') = 0$ for all integers $r' > r$.  
\ec 

\begin{proof} 
This follows immediately from Theorem~\ref{PhaseTransition:theorem:fundamental}.
\end{proof}

Suppose that  $A$ is an infinite generating set for a group $G$ such that  the word growth function satisfies $L_A(r') = 0$ for some integer $r'$.  Let $r$ be the largest integer such that $L_A(r) > 0$ and let $s = L_A(r) = |S_A(r)|$.  The ordered pair $(r,s)\in \N \times (\N \cup \{\infty\})$ is called the \emph{phase transition} of the generating set $A$, and $S_A(r)$ is called the \emph{transition set}.   For example, let $G = \Z^r$ be the additive group of $r$-dimensional lattice points.  Let $A_i$ be the set of all lattice points with an integer in the $i$th place and 0's elsewhere, and let $A = \bigcup_{i=1}^r A_i$.  The word length of the lattice point $x\in \Z^r$ is the number of nonzero coordinates of $x$.  It follows that $A$ is an infinite generating set for $\Z^r$ with phase transition $(r,\infty)$ and transition set 
\[
\{ (x_1,\ldots, x_r) \in \Z^r : x_i \neq 0 \text{ for all } i=1,\ldots, r\}.
\]
We shall construct, for every $(r,s) \in \N \times \N$, a group $G$ and an infinite generating set $A$ for $G$ such that $A$ has phase transition $(r,s)$.

\section{Direct products of groups}

We define $0\cdot\infty = \infty\cdot 0 = 0$ and $k\cdot\infty = \infty\cdot k = \infty$ for $k \in \N \cup \{ \infty\}$.   
This is consistent with the fact that the direct product of an infinite set with the empty set is the empty set and the direct product of an infinite set with a nonempty set is an infinite set.

\bl     \label{PhaseTransition:lemma:DirectProduct}
Let $G_1$ and $G_2$ be groups with identity elements $e_1$ and $e_2$ and nonempty generating sets $A_1$ and $A_2$, respectively.  The sets $A_1$ and $A_2$ can be finite or infinite.  Let $G = G_1 \times G_2$ and let
\[
A = \left(A_1 \times \{ e_2 \}  \right) \cup \left( \{ e_1 \} \times A_2 \right).
\]
The set $A$ generates $G$, and 
\[
\ell_A(x_1,x_2) = \ell_{A_1}(x_1) + \ell_{A_2}(x_2)
\]
for all $(x_1,x_2) \in G_1 \times G_2$.
Moreover, 
\[
S_A(r) = \bigcup_{r'=0}^r \left( S_{A_1}(r') \times S_{A_2}(r-r') \right) 
\]
and
\beq    \label{PhaseTransition:DirectProductGrowth}
L_A(r) = \sum_{r'=0}^r L_{A_1}(r') L_{A_2}(r-r').
\eeq
\el

\begin{proof}
We observe that the sets $A_1 \times \{ e_2 \}$ and $ \{ e_1 \} \times A_2 $ commute with each other, and so every product of elements of $A$ can be written as a product in which the elements of $A_1 \times \{ e_2 \}$ precede, that is, occur to the left of, the elements of $ \{ e_1 \} \times A_2 $.

Suppose that $x_1 \in G_1$ and $x_2 \in G_2$ satisfy $\ell_{A_1}(x_1) = r_1$ and $\ell_{A_2}(x_2) = r_2$.  Choose  $a_{i,1} \in A_1$ for $i=1,\ldots, r_1$ such that $x_1 = a_{1,1}a_{2,1}\cdots a_{r_1,1}$, and choose  $a_{i,2} \in A_2$ for $i=1,\ldots, r_2$ such that $x_2 = a_{2,1}a_{2,2}\cdots a_{r_2,2}$.  Then
\begin{align*}
(x_1,x_2) 
& = ( a_{1,1}a_{2,1}\cdots a_{r_1,1}, a_{2,1}a_{2,2}\cdots a_{r_2,2})\\
& =  (a_{1,1},e_2)(a_{2,1},e_2)\cdots (a_{r_1,1},e_2)(e_1,a_{2,1})(e_1,a_{2,2})\cdots (e_1,a_{r_2,2})
\end{align*}
and so $\ell_A(x_1,x_2) \leq r_1 + r_2 = \ell_{A_1}(x_1) + \ell_{A_2}(x_2)$.

Conversely, suppose that there exist $(a_{i,1},e_2) \in A$ for $i=1,\ldots, s_1$ and $(e_1,a_{i,2}) \in A$ for $i =1,\ldots, s_2$ 
such that 
\[
(x_1,x_2)  =  (a_{1,1},e_2)(a_{2,1},e_2)\cdots (a_{s_1,1},e_2)(e_1,a_{2,1})(e_1,a_{2,2})\cdots (e_1,a_{s_2,2}).
\]
The length of this word is $s_1+s_2$.  We have
\[
(x_1,x_2) 
 =  (a_{1,1}a_{2,1} \cdots a_{s_1,1},a_{2,1}\cdots a_{2,2}a_{s_2,2})
\]
and so $x_1 = a_{1,1}a_{2,1} \cdots a_{s_1,1}$ and $x_2 = a_{2,1}a_{2,2} \cdots a_{s_2,2}$.  It follows that $s_1 \geq \ell_{A_1}(x_1)$ and $s_2 \geq \ell_{A_2}(x_2)$, hence $\ell_A(x_1,x_2) \geq \ell_{A_1}(x_1) + \ell_{A_2}(x_2)$.  Therefore, $\ell_A(x_1,x_2) = \ell_{A_1}(x_1) + \ell_{A_2}(x_2)$.  

This identity implies that $(x_1,x_2) \in S_A(r)$ if and only if $\ell_{A_1}(x_1) = r'$ and $\ell_{A_2}(x_2) = r - r'$ for some $r' \in \{0,1,2,\ldots, r\}$.  Equivalently, $(x_1,x_2) \in S_A(r)$ if and only if $x_1 \in S_{A_1}(r')$ and $x_2 \in S_{A_2}(r-r')$ and so 
$(x_1,x_2) \in S_{A_1}(r') \times S_{A_2}(r-r')$  for some $r' \in \{0,1,2,\ldots, r\}$, that is, 
\[
 S_A(r) = \bigcup_{r'=0}^r \left( S_{A_1}(r') \times S_{A_2}(r-r') \right).
\]
Since the sets $S_{A_1}(r') \times S_{A_2}(r-r')$ are pairwise disjoint for $r' \in \{0,1,2,\ldots, r\}$, it follows that 
\begin{align*}
L_A(r) 
& = \left| \bigcup_{r'=0}^r \left( S_{A_1}(r') \times S_{A_2}(r-r') \right)\right| \\
& =\sum_{r'=0}^r  \left|  S_{A_1}(r') \times S_{A_2}(r-r')  \right| \\
& =\sum_{r'=0}^r  \left|  S_{A_1}(r') \right| \left|S_{A_2}(r-r')  \right| \\
& =\sum_{r'=0}^r  L_{A_1}(r') L_{A_2}(r-r').
\end{align*}
This completes the proof.  
\end{proof}

\bt     \label{PhaseTransition:theorem:DirectProduct}
Let $G_1$ and $G_2$ be groups with infinite generating sets $A_1$ and $A_2$, and with  phase transitions $(r_1,s_1)$ and $(r_2,s_2)$, respectively.   Consider the generating set  $A = (A_1\times \{e_2\}) \cup (\{e_1\} \times A_2)$ for the group $G_1 \times G_2$.    The phase transition of  $A$ is $(r_1+r_2, s_ 1s_ 2)$ and the transition set is $S_{A_1}(r_1) \times S_{A_2}(r_2)$. 
\et

\begin{proof}
By Lemma~\ref{PhaseTransition:lemma:DirectProduct}, the growth function $L_A$ satisfies the polynomial recurrence~\eqref{PhaseTransition:DirectProductGrowth}.
Since $0 \cdot \infty = \infty \cdot 0 = 0$, we have
$
L_{A_1}(r') L_{A_2}(r-r') = 0
$
if and only if $r' > r_1$ or $r-r' > r_2$. 
If $r > r_1+r_2$ and $r' \leq r_1,$ then $r-r' \geq r-r_1 > r_2$.  This implies that $L_A(r) = 0$ if $r>r_1+r_2$.  
Let $r = r_1 + r_2$.  Since  
$
L_{A_1}(r') L_{A_2}(r-r') = 0
$
unless $r' = r_1$, it follows that 
\begin{align*}
L_A(r) & = \sum_{r'=0}^{r_1+r_2} L_{A_1}(r') L_{A_2}(r_1+r_2-r') \\
& = L_{A_1}(r_1) L_{A_2}(r_2) \\
& = s_ 1s_ 2 > 0.
\end{align*}
We apply Theorem~\ref{PhaseTransition:theorem:fundamental} to complete the proof.
\end{proof}

\bl      \label{PhaseTransition:lemma:DirectProduct2}
Let $G_1$ and $G_2$ be groups with  identity elements $e_1$ and $e_2$ and nonempty generating sets $A_1$ and $A_2$, respectively.  
The sets $A_1$ and $A_2$ can be finite or infinite. 
Let $G = G_1 \times G_2$ and let
\[
A = (A_1 \times A_2) \cup \left(A_1 \times \{ e_2 \}  \right) \cup \left( \{ e_1 \} \times A_2 \right).
\]
The set $A$ generates $G$, and 
\[
\ell_A(x_1,x_2) = \max\left( \ell_{A_1}(x_1) , \ell_{A_2}(x_2) \right)
\]
for all $(x_1,x_2) \in G_1 \times G_2$.  
Moreover, 
\[
S_A(r) = \left( S_{A_1}(r)\times S_{A_2}(r) \right)
 \cup \bigcup_{r'=0}^{r-1} \left( S_{A_1}(r')\times S_{A_2}(r)\right)  \cup  \bigcup_{r'=0}^{r-1} \left( S_{A_1}(r)\times S_{A_2}(r') \right) \]
and 
\beq    \label{PhaseTransition:DirectProductGrowth2}
L_A(r) = L_{A_1}(r)L_{A_2}(r) + L_{A_2}(r) \sum_{r'=0}^{r-1} L_{A_1}(r')  +  L_{A_1}(r)  \sum_{r'=0}^{r-1} L_{A_2}(r').
\eeq
\el

\begin{proof}
Suppose that $x_1 \in G_1$ and $x_2 \in G_2$ satisfy $\ell_{A_1}(x_1) = r_1$ and $\ell_{A_2}(x_2) = r_2$.  Choose  $a_{i,1} \in A_1$ for $i=1,\ldots, r_1$ such that $x_1 = a_{1,1}a_{2,1}\cdots a_{r_1,1}$, and choose  $a_{i,2} \in A_2$ for $i=1,\ldots, r_2$ such that $x_2 = a_{2,1}a_{2,2}\cdots a_{r_2,2}$.  If $r_1 \leq r_2$, then
\begin{align*}
(x_1,x_2) 
& = ( a_{1,1}a_{2,1}\cdots a_{r_1,1}, a_{2,1}a_{2,2}\cdots a_{r_2,2})\\
& =  (a_{1,1},a_{2,1}) (a_{2,1},a_{2,2}) \cdots (a_{r_1,1},a_{r_1,2})(e_1,a_{r_1+1,2})(e_1,a_{r_1+2,2})\cdots (e_1,a_{r_2,2})
\end{align*}
and so $\ell_A(x_1,x_2) \leq r_2 = \max\left(\ell_{A_1}(x_1), \ell_{A_2}(x_2)\right)$.

Conversely, suppose that there exist $(a_{i,1},a_{i,2}) \in A$ for $i =1,\ldots, r$  such that 
\[
(x_1,x_2)  =  (a_{1,1},a_{2,1})(a_{2,1},a_{2,2})\cdots (a_{r,1},a_{r,2}).
\]
The length of this word is $r$.  
We have
\[
(x_1,x_2)   =  (a_{1,1}a_{2,1} \cdots a_{r,1},a_{2,1} a_{2,2} \cdots a_{r,2})
\]
and so $x_1 = a_{1,1}a_{2,1} \cdots a_{r,1}$ and $x_2 = a_{2,1}a_{2,2} \cdots a_{r,2}$.   It follows that 
\[
r \geq \max\left( \ell_{A_1}(x_1), \ell_{A_2}(x_2)\right)
\]
and so $\ell_A(x_1,x_2) \geq \max\left( \ell_{A_1}(x_1), \ell_{A_2}(x_2)\right)$.  Therefore, 
\[
\ell_A(x_1,x_2) = \max\left( \ell_{A_1}(x_1), \ell_{A_2}(x_2)\right).
\]
The rest of the proof is similar to the proof of Lemma~\ref{PhaseTransition:lemma:DirectProduct}.
\end{proof}

\bt    \label{PhaseTransition:theorem:DirectProduct2}
Let $G_1$ and $G_2$ be groups with infinite generating sets $A_1$ and $A_2$, and with  phase transitions $(r_1,s_ 1)$ and $(r_2,s_ 2)$, respectively.   Consider the generating set  $A = (A_1 \times A_2) \cup (A_1\times \{e_2\}) \cup (\{e_1\} \times A_2)$ for the group $G_1 \times G_2$.    
The phase transition of $A$ is $(\max(r_1,r_2), \infty)$. 
\et

\begin{proof}
Let $r_1 \leq r_2$.  If $r > r_2$, then $L_{A_1}(r) = L_{A_2}(r) = 0$ and 
\[
L_A(r) = L_{A_1}(r)L_{A_2}(r) + L_{A_2}(r) \sum_{r'=0}^{r-1} L_{A_1}(r')  +  L_{A_1}(r)  \sum_{r'=0}^{r-1} L_{A_2}(r') = 0.
\]
We have
\[
L_A(r_2) \geq L_{A_1}(r_1)L_{A_2}(r_2) = s_ 1s_ 2 > 0.
\]
If $s_ 1 < \infty$ and $s_ 2 < \infty$, then $r_2 \geq r_1 \geq 2$. 
It follows from Theorem~\ref{PhaseTransition:theorem:fundamental} that $L_{A_1}(r_1-1) = \infty$, and so 
\[
L_A(r_2) \geq L_{A_1}(r_1-1)L_{A_2}(r_2) = \infty.
\]
Therefore, $(\max(r_1,r_2),\infty)$ is the phase transition of the generating set $A$.  This completes the proof.
\end{proof}

\section{Constructions of groups and generating sets with arbitrary phase transitions}

In this section we construct, for every ordered pair $(r,s)$ of positive integers with $r \geq 2$, a group $G$ and an infinite generating set $A$  whose phase transition is $(r,s)$.  For even integers $s$, we construct the appropriate generating sets for the additive group \Z.  Constructions of transition sets for odd integers $s$ are more complicated.  Recall that if the set $A$ generates the group $G$, then $\ell_A(x) = \ell_A(x^{-1})$ for all $x\in G$.  It follows that a transition set with an odd number of elements must contain an element $x \neq e$ such that $x^{-1} = x$, that is, an element of order 2, also called an \emph{involution}.  The group of integers is torsion-free, and, in particular, contains no involutions.    However,  the group $\Z \times (\Z/2\Z)$  does contain an involution; the element $(0,\overline{1})$ is the unique involution in this group.  
For every integer $r \geq 2$ and every odd positive integer $s$, we construct a generating set on the finite direct product  $\left( \Z \times (\Z/2\Z)\right)^{[r/2]}$ with phase transition $(r,s)$.

\bt                \label{PhaseTransition:theorem:Zphases}
Let $r$ and $s$ be positive integers such that $r \geq 2$ and $s$ is even.  Let $W$ be a set of $s/2$ positive integers.  There  exists an infinite set $A$ of positive integers such that $A$ generates \Z\ and $S_A(r) = W \cup (-W)$.   In particular, the pair $(\Z,A)$ has phase transition $(r,s)$.
\et

\begin{proof}
If $r = 2$, then the set $A = \N \setminus W$ is the unique set of positive integers such that $S_A(2) = W \cup (-W)$.  The pair $(\Z,A)$ has phase transition $(2,s)$.

Let $r \geq 3$ and $a_0 = 0$.    Let $\{a_{2i-1}\}_{i=1}^{\infty}$ be  a  sequence of positive integers such that $a_1 > \max(W)$ and 
\[
a_{2i+1} >  (r-2)  \left( (r-2) a_{2i-1} + i \right) + \max(W)
\]
for all $i \geq 1$.  We define
\[
a_{2i} = (r-2) a_{2i-1} + i.
\] 
Then 
\[
a_{2i-1} < a_{2i} <  (r-2)a_{2i} + \max(W) <  a_{2i+1}
\]
for all $i \in \N$.  Consider the set 
\[
A = \bigcup_{\substack{i=1\\i\notin W}}^{\infty} \{a_{2i-1}, a_{2i} \}.
\]
For all $i \in \N \setminus W$ we have  
\[
i  = a_{2i} - (r -2)a_{2i-1}  \in (r-1) (A \cup (-A))
\]
and so $\ell_A(i) \leq r-1$ for all $i \notin  W \cup (-W)$.  Equivalently,
\[
 \Z \setminus (W \cup (-W)) \subseteq \bigcup_{r'=0}^{r-1} S_A(r').
\]
We shall prove that 
\[
\Z \setminus (W \cup (-W)) = \bigcup_{r'=0}^{r-1} S_A(r').
\]
Let $n$ be a positive integer such that $n \in \bigcup_{r'=0}^{r-1} S_A(r')$.  There exist finite disjoint multisets $J$ and $K$ of positive integers such that 
\[
n = \sum_{j\in J}a_j  - \sum_{k\in K}a_k
\]
and $|J|+|K| \leq r-1$, where $|J|$ and $|K|$ denote the cardinalities of the multisets $J$ and $K$.  Since $n$ is positive and $a_{i+1} > (r-2)a_i$ for all $i$, it follows that $\max(J \cup K) \in J$.   
If $\max(J \cup K) = 2i+1$, then 
\[
n  \geq a_{2i+1} - \sum_{k\in K}a_k 
\geq a_{2i+1} - (r-2)a_{2i}  > \max(W)
\]
and so $n \notin W$.  

Suppose that $\max(J \cup K) = 2i$.
If $J = \{2i\}$ and if $K$ is the multiset that contains only the integer $a_{2i-1}$  with multiplicity $r-2$, then 
\[
n =  a_{2i} - \sum_{k\in K}a_k =a_{2i} - (r-2)a_{2i-1} = i \notin W.  
\]
Otherwise,
\begin{align*}
n & \geq  a_{2i} - \sum_{k\in K}a_k  \geq a_{2i} - (r-3)a_{2i-1} - a_{2i-2} \\
& = a_{2i-1} - a_{2i-2} + i.
\end{align*}
If $i=1$, then $n > a_1 > \max(W)$.
If $i \geq 2$, then 
\[
n > a_{2i-1} - a_{2i-2} = (r-3) a_{2i-2} +  \max(W)  \geq  \max(W).
\]
Therefore,  $n \notin W$ and 
\[
\Z \setminus \bigcup_{r=0}^{r-1} S_A(r) = W \cup (-W).
\]
Since this set is finite and nonempty, Theorem~\ref{PhaseTransition:theorem:fundamental} implies that 
\[
S_A(r) = W \cup (-W)
\]
and $A$ has phase transition $(r,s)$.
This completes the proof.
\end{proof}

In the cyclic group $\Z/2\Z$ we denote the residue classes $2\Z$ and $1+2\Z$ by $\overline{0}$ and $\overline{1}$, respectively.

\bl   \label{PhaseTransition:lemma:phase transition-2-s}
Consider the additive abelian group $\Z \times (\Z/2\Z)$.  
For every positive odd integer $s $, let $t = (s+1)/2$ and 
\[
A_s = \{ (k,\overline{0})  : k = 1,2,3,\ldots\} \cup \{ (k,\overline{1}) : k = t,t+1,t+2,\ldots\}.
\]
Then $A_s$ generates $\Z \times (\Z/2\Z)$ and has phase transition $(2,s)$.
\el

\begin{proof}
It suffices to observe that 
\[
S_{A_s}(0) \cup S_{A_s}(1) 
= \Z \times (\Z/2\Z) \setminus \{ (k, \overline{1}) : k = 0,\pm 1, \ldots , \pm (t-1) \}.
\]
This completes the proof.
\end{proof}

\bl     \label{PhaseTransition:lemma:phase transition-3-s}
Consider the additive abelian group $\Z \times (\Z/2\Z)$.    
Let $s$ be a positive odd integer, and $t = (s+1)/2$.
Let
$\{a_k\}_{k=1}^{\infty}$ be a sequence of positive integers such that 
\[
a_k + 2k + t < a_{k+1}
\]
for all $k \geq 1$.  
For every positive integer $s$, let
\[
A_s = \{ (a_k,\overline{0}) , (a_k + k,\overline{0})  : k = 1,2,3,\ldots\} \cup \{ (a_k+2k,\overline{1}) :  k = t,t+1,t+2,\ldots \}.
\]
Then $A_s$ generates $\Z \times (\Z/2\Z)$ and has phase transition $(3,s)$.
\el

\begin{proof}
For all positive integers $k$ we have 
\[
(k,\overline{0}) = (a_k + k,\overline{0}) - (a_k,\overline{0})
\]
and so $\ell_A( (\pm k,\overline{0})) \leq 2$.  On the other hand, $\ell_A(a_k+k+i, \overline{0}) \geq 2$ for all $k \geq 2$ and $i=1,2,\ldots, k-1$, and $L_{A_s}(2) = \infty$.

For every integer $k \geq t$ we have 
\[
(k,\overline{1}) = (a_k + 2k,\overline{1}) - (a_k+k,\overline{0})
\]
and so $\ell_A(( \pm k,\overline{1})) \leq 2$.  Theorem~\ref{PhaseTransition:theorem:fundamental} implies that 
\begin{align*}
S_{A_s} (3) 
& = \Z \times (\Z/2\Z) \setminus \left( S_{A_s} (0) \cup S_{A_s} (1) \cup S_{A_s} (2) \right) \\
& \subseteq \{ (\pm k,\overline{1}) : k = 0,1,2,\ldots, t-1 \}.
\end{align*}
We shall prove that $S_{A_s} (3) = \{ (\pm k,\overline{1}) : k = 0,1,2,\ldots, t-1 \}$.  It suffices to show that the sets $S_{A_s} (2)$ and 
$\{ (k,\overline{1}) : k = 0,1,2,\ldots, t-1 \}$ are disjoint.

If not, then there exists $k \in \{0,1,2,\ldots, t-1 \}$  such that $(k,\overline{1}) \in S_{A_s}(2)$.  This means that there exist positive integers $i$ and $j$ with $j \geq t$ and $\varepsilon \in \{ 0,1 \}$ such that
\[
(k,\overline{1}) = \pm (a_j+2j,\overline{1}) \pm (a_i+\varepsilon i,\overline{0})
 = (\pm(a_j+2j)  \pm(a_i+\varepsilon i),\overline{1})
 \]
 and so
 \[
 \pm(a_j+2j)  \pm(a_i+\varepsilon i) =k.
 \]
 
Since $a_j+2j + a_i+\varepsilon i > 2j \geq 2t > k$, it follows that either 
(i) $k = a_j+2j - ( a_i+\varepsilon i )$ or 
(ii) $k = a_i+\varepsilon i - ( a_j+2j )$.  
In case (i), we have $j \geq i$ and 
\[
k = a_j+2j - ( a_i+\varepsilon i ) \geq a_j+2j - ( a_j+j ) = j \geq t > k
\]
which is absurd.  
In case (ii) we have $i > j$ and 
\[
k = a_i+\varepsilon i - ( a_j+2j ) \geq a_i - ( a_{i-1}+2(i-1) ) > t \geq k
\]
which is also absurd.  It follows that 
\[
S_{A_s} (3) = \{ (\pm k,\overline{1}) : k = 0,1,2,\ldots, t-1 \}
\]
and so the group $\Z \times (\Z/2\Z)$ with infinite generating set $A_s$  has phase transition $(3,s)$.
\end{proof}

\bt   \label{PhaseTransition:theorem:OddPhase}
For positive integers $r$ and $s$ with $r \geq 2$ and $s$ odd, the additive  group $\left(  \Z \times (\Z/2\Z)  \right)^{[r/2]}$ has an infinite  generating set  $A$ with phase transition $(r,s)$.
\et

\begin{proof}
In Lemma~\ref{PhaseTransition:lemma:phase transition-2-s}, for every odd positive integer $s$ we constructed an infinite generating set $A_s$ for the group $\Z \times (\Z/2\Z)$ with phase transition $(2,s)$.  Applying this in the case $s=1$ to Theorem~\ref{PhaseTransition:theorem:DirectProduct}, we see that if a group $G$ has phase transition $(r,s)$ with respect to some infinite generating set, then there also exists an infinite generating set  for the group $G\times \Z \times (\Z/2\Z) $ with phase transition $(r+2,s)$.  If we start with the generating set $A_s$ for the group $\Z \times (\Z/2\Z)$ with phase transition $(2,s)$, then for every even positive integer $r$ we obtain inductively an infinite generating set for the group  $\left(  \Z \times (\Z/2\Z)  \right)^{r/2}$ with phase transition $(r,s)$.  

Similarly, we obtain from Lemma~\ref{PhaseTransition:lemma:phase transition-3-s} 
a generating set $A_s$ for the group $\Z \times (\Z/2\Z)$ with phase transition $(3,s)$.  Again applying  Theorem~\ref{PhaseTransition:theorem:DirectProduct},  for every odd positive integer $s$ we obtain inductively an infinite generating set for the group $\left(  \Z \times (\Z/2\Z)  \right)^{[r/2]}$ with phase transition $(r,s)$.    This completes the proof.
\end{proof}

\section{Open Problems}

\bprob
Let $G$ be a group.  Classify the set of all finite phase transitions and finite transition sets associated with infinite generating sets for $G$.  
\eprob

Theorem~\ref{PhaseTransition:theorem:Zphases} solves this problem for $G = \Z$:  The set of finite phase transitions is $\{(r,s) \in \N^2 : r \geq 2 \text{ and $s$ even} \}$, and for every $r \geq 2$, the set of finite transition sets of integers is 
\[
 \{ W \cup (-W) : W \subseteq \N \text{ and } 0 < |W| < \infty \}.
\]
However, the problem of describing the finite transition sets in an arbitrary group, even an abelian group, is difficult.  For example, consider the group  
\[
G = \bigoplus_{i=0}^{\infty} G_i
\]
where $G_0 = \Z/4\Z$ and $G_i = \Z/2\Z$ for $i \in \N$.  Let 
$w^* = (w_0,w_1,w_2,\ldots)$, where $w_0 = 2 + 4\Z$ and $w_i=2\Z$ for $i\in \N$.  Let $W$ be any finite symmetric subset of $G\setminus \{e\}$ with $w^* \in W$.   We shall prove that, for every $r \geq 3$,  there does not exist an infinite set $A$ that generates $G$ and has the associated  transition set $S_A(r) = W$.  If there  did exist such a generating set $A$, then $\ell_A(w^*) \geq 3$.  
Let $x = (x_0,x_1,x_2,\ldots)$ be an element of $G\setminus W$ with  $x_0 = 1 + 4\Z$.  Then $1 \leq \ell_A(x) \leq r-1$, and there is a sequence of at most $r-1$ elements of $A \cup (-A)$ whose sum is $x$.  Since $x_0 = 1 + 4\Z$, at least one of these summands must be of the form $a = (a_0,a_1,a_2,\ldots)$ with $a _0 = 1 + 4\Z \text{ or } 3 + 4\Z$, and so $2a = w^* \in S_A(2)$, that is, $\ell_A(w^*) \leq 2$.  This is a contradiction.

\bprob
This is the \emph{inverse problem for finite phase transitions}. 
Let $G$ be a group, and let $(r,s)$ be a finite phase transition for $G$ with transition set $S$.  Describe the set of all infinite generating sets $A$ for $G$ such that $A$ has  transition set $S$.
\eprob

\bprob
The \emph{counting function} $A(t)$ of a set $A$ of integers counts the number of positive elements of $A$ that do not exceed $t$.  
For $r \geq 2$, let $W$ be a nonempty finite  set of positive integers and let $A$ be a set of positive integers such that $A$ generates \Z\ and $S_A(r) = W \cup (-W)$.  
Compute the greatest density of a set $A$ whose transition set is $W \cup (-W)$.  In particular, estimate or compute 
\[
\sup\left(\{\theta > 0: \text{there exists $A \subseteq \N$ such that $S_A(r) = W\cup (-W)$ and $A(t) \gg t^{\theta}$ } \} \right).
\]
\eprob

\bprob
Determine the set of phase transitions and transition sets for the group $\Z \times (\Z/2\Z)$.
\eprob

\bprob
Let $G$ be a group.  Classify the set of all infinite phase transitions and infinite transition sets associated with infinite generating sets for $G$.  This is not known even for $G = \Z$.
\eprob

\bprob
This is the \emph{inverse problem for infinite phase transitions}. 
Let $G$ be a group, and let $(r,\infty)$ be an infinite phase transition for $G$ with transition set $S$.  Classify the set of all infinite generating sets $A$ for $G$ such that $A$ has phase transition $(r,\infty)$  and transition set $S$.
\eprob

\bprob
Let $G$ be an infinite group and let $\mathcal{A}$ be the set of all infinite generating sets for $G$.  Is there an algorithm to determine if a generating set $A \in \mathcal{A}$ has a phase transition?
\eprob

\bprob
Let $G$ be an infinite group.   Is there a method to determine if an infinite generating set $A$ for $G$ has a finite or an infinite phase transition?
\eprob

\emph{Acknowledgement.} I thank Jacob Fox for many helpful discussions about this work.

\end{document}